\documentclass[12pt]{article}
\usepackage{amsmath,amsfonts,amssymb,amscd}
\usepackage[english]{babel}
\newcommand{\PP}{{\mathbb{P}}}

\newcommand{\C}{{\mathbb{C}}}
\newcommand{\Q}{{\mathbb{Q}}}

\newcommand{\calO}{{\cal O}}

\newcommand{\calM}{{\cal M}}
\newtheorem{theorem}{Theorem}

\newtheorem{lemma}[theorem]{Lemma}

\newcommand{\iso}{\stackrel{\sim}{\rightarrow}}
\newcommand{\into}{\hookrightarrow}
\title{Degeneration of the Leray spectral sequence for certain geometric
quotients}
\author{C.A.M.~Peters\\Department of Mathematics, University of Grenoble I\\UMR 5582 CNRS-UJF, 38402-Saint-Martin d'H\`eres\\France, email: peters@ujf-grenoble.fr \and J.H.M.~Steenbrink\thanks{The second author thanks the University of Grenoble I for its hospitality and financial support}\\Department of Mathematics, University of Nijmegen\\
Toernooiveld, NL-6525 ED Nijmegen\\The Netherlands, email: steenbri@sci.kun.nl}

\begin{document}
\maketitle

\begin{abstract}
We prove that the Leray spectral sequence in rational cohomology for the quotient map $U_{n,d} \to U_{n,d}/G$ where $U_{n,d}$ is the affine variety of equations for smooth hypersurfaces of degree $d$ in $\PP^n(\C)$ and $G$ is the general linear group, degenerates at $E_2$. 
\end{abstract}
\section{Introduction}

We consider an affine complex algebraic group $G$ which acts on a smooth
algebraic variety $X$. Assume that a geometric quotient $f:X\to Y$ for 
the action of $G$ on $X$ exists (cf. \cite[Sect. 0.1]{MFK}).  We want to give
geometric conditions ensuring that the Leray spectral sequence degenerates at
$E_2^{p,q}=H^p(Y,R^qf_*\Q)$.

The cohomology ring of $G$ is well known (\cite{Ho}, \cite{B}). It is an exterior
algebra with exactly one generator $\eta_i$ in certain odd degrees $2r_i-1$,
$i=1,\dots, r=r(G)$, the {\it rank} of $G$.  So, if the Leray spectral sequence
degenerates at $E_2$, knowing the cohomology of the source $X$ is equivalent to
knowing that of the target $Y$.  As an example of how this could be used, we point out
that for any group $G$ acting with finite stabilizers on a topological space $X$
the equivariant cohomology $H^\bullet_G(X,\Q)$ equals $H^\bullet(X/G,\Q)$ (\cite[\S1, Remark 2]{Br}) and the former can often be calculated group theoretically. See
\cite{Br} for examples. So, in these cases one knows $H^\bullet(X,\Q)$.

We prove a general result (Theorem \ref{maintheorem}) giving sufficient
geometric conditions for this to happen. These turn out to be satisfied for the
group  $\mathrm{GL}_{n+1}(\C)$ acting  on the affine variety $U_{n,d}$ of
those homogeneous polynomials of degree $d$ in $(n+1)$ variables which give smooth
hypersurfaces in $\PP^n$:

\begin{theorem} 
Let $d \geq 3$. Then the Leray spectral sequence in  rational cohomology for
the quotient  map $U_{n,d} \to U_{n,d}/G$, $G=\mathrm{GL}_{n+1}(\C)$ degenerates
at $E_2$.
\end{theorem}

\vspace{4mm}\noindent {\bf Examples} 
\begin{enumerate}
\item By results of Vassiliev \cite{V}  the map $H^\ast(U_{n,d};\Q) \to
H^\ast(\mathrm{GL}_{n+1}(\C);\Q)$ is an isomorphism in the cases
$(n,d)=(2,3),(3,3)$. Moreover Gorinov \cite{G} has proved the same
result for the cases $(n,d)=(4,3),(2,5)$. It follows that $M_{n,d}$ has the
rational cohomology of a point in these cases. 
\item For the case $(n,d)=(2,4)$  it follows from \cite{V} and Theorem \ref{LH}
that the space $M_{2,4}$ has a cohomology group of dimension $1$ in degrees $0$
and $6$ and has zero rational cohomology in other degrees.  This agrees with a
result of Looijenga \cite{L} about the Poincar\'e-Serre polynomial of
$M_{2,4}$: 
$$H^6(M_{2,4};\Q) \simeq \Q(-6)$$ 
and the other cohomology groups are those of a
point. 
\end{enumerate}

\noindent{\bf Remark}
In \cite{ACT} there is a description of $M_{3,3}$ using periods of
threefolds. This moduli space turns out to be a certain explicitly described open 
subset of the quotient of complex hyperbolic $4$-space by a certain discrete
group.   From this description it is quite unexpected that $M_{3,3}$ has the
rational cohomology of a point. It is an interesting question to calculate the
cohomology of the various compactifications of $M_{3,3}$ studied in loc. cit.

\section{Generalizing the Leray-Hirsch theorem}   \label{one}

The proof of the  Leray-Hirsch theorem  as given on \cite[ p. 229]{H} is valid
for a locally trivial fibration  $p:M\to B$. For cohomology with {\it rational}
coefficients, the same proof applies to a slightly more general situation: 

\hskip1ex

\noindent{\bf Definition}  A continuous map $p:M\to B$ is a locally trivial
fibration, say with fibre $F$, in the {\it orbifold sense\/} if for every $b\in
B$  there exists a neigbourhood $V_b$, a topological space $U_b$, and a
topological group
$G_b$ such that

\begin{enumerate} 

\item  $G_b$ acts on $U_b$ and on $F$; the action on $F$ is
by homeomorphisms homotopic to the identity;

\item  $V_b$ is homeomorphic to $U_b/G_b$;

\item   $p^{-1}V_b$ is homeomorphic to the quotient of $U_b\times F$ by the
product action of $G_b$.
\end{enumerate}

In this setting, composing the natural quotient map $F\to F/G_b$  with the
homeomorphism  $(F/G_b)\iso  p^{-1}b$ and the inclusion $p^{-1}b\into X$, defines the {\it orbifold fibre inclusion} $r_b:F\to X$.

\hskip1ex 

Indeed, in this setting  the proof as
given in loc. cit. applies starting from the observation that over the rationals
we still have graded isomorphisms (replacement of the K\"unneth formula) 
\begin{eqnarray*}
H^\bullet(p^{-1}V_b;\Q)  & \cong H^\bullet(U_b\times F;\Q)^{G_b} &\cong
 H^\bullet(U_b;\Q)^{G_b}\otimes H^\bullet(F;\Q)^{G_b}   \\
&& \cong H^\bullet(V_b;\Q)\otimes H^\bullet(F;\Q),
\end{eqnarray*}
because $g\in G_b$  acts trivially on $H^q(F;\Q)$ since it is
homotopic to the identity by assumption.

We thus arrive at: 

\begin{theorem} \label{LH} Let $p:M\to B$ be a   fibration which is locally
trivial in the orbifold sense. Suppose that for all $q\ge 0$ there exist classes
$e_1^{(q)},\dots,e_{n(q)}^{(q)}\in H^q(M;\Q)$ that restrict to a basis for
$H^q(F;\Q)$ under the map induced by the orbifold fibre inclusion $r_b:F \to  M$
. 
The map  $a\otimes r_b^*(e_i)\mapsto p^*a\cup e_i$, $a\in H^\bullet(B;\Q)$ 
extends linearly to a graded linear isomorphism
$$
H^\bullet (B;\Q)\otimes H^\bullet (F;\Q)\iso H^\bullet (M;\Q).
$$
\end{theorem}

\vspace{4mm}\noindent {\bf Example} 
Let $\phi:X \to Y$ be a  geometric quotient for $G$. Suppose that $G$ is
connected and that for all $x\in X$, the identity component of the stabiliser
$S_x$ of $x$ is contractible (e.g. when $S_x$ is finite). 
For  $y\in Y$ we take for $U_y$ any open slice for the action of $G$ through $x\in
\phi^{-1}y$, i.e. a  contractible submanifold  through $x$ which intersects
$Gx$ transversally at $x$. Then, if $gx$ is any other point in same orbit,
$gU_y$ is a slice through $gx$ and $gS_xg^{-1}=S_{gx}$ so that for all
$g\in G$, the quotient $gU_{y}/S_{gx}$ gives the same neighbourhood $V_y$ of $y$
.
We have $(U_{y}\times G)/S_x=\phi^{-1}(V_y)$. The assumption that $G$
is connected implies that multiplication by $g\in G$ is homotopic to the identity
in $G$. So $\phi$ is indeed locally trivial in the orbifold sense (with typical
fibre $G$).
 
We  study this example in more detail in the next section.

\section{The case of a geometric quotient for a reductive group}

We  assume that  $G$ is a reductive complex affine group, that $V$ is a
representation space for $G$ and that $X$ is an open set of stable points with
complement $\Sigma=V\setminus X$.  For $x\in X$  the orbit map is denoted as
follows
\begin{eqnarray*}
o_x:G &\to X \\
              g&\mapsto g(x), 
\end{eqnarray*}
and the geometric quotient by
$$
\phi: X\longrightarrow  Y=X/G.
$$
Recall that $H^\bullet(G)$ is an exterior algebra freely generated by classes 
$\eta_i\in H^{2r_i-1}(G)$. Note also that $V$ being a vector space, we have
isomorphisms
$$
H^{2r_i-1}(X)\iso H^{2r_i}_{\Sigma}(V).
$$
We can now apply the variant of the Leray-Hirsch theorem as stated in the previous
section to the geometric quotient $\phi$  and we obtain:

\begin{theorem} \label{maintheorem}
Suppose that there are schemes $Y_i\subset \Sigma$ of pure codimension $r_i$ in
$V$ whose fundamental classes map to a non-zero multiple of $\eta_i$ under the
composition
$$
H^{2r_i}_{Y_i}(V) \to  H^{2r_i}_{\Sigma}(V) \iso H^{2r_i-1}(X)
{\stackrel{o_x^*}{\longrightarrow}} H^{2r_i-1}(G).
$$
Denote the image of $[Y_i]$ in $H^\bullet(X;\Q)$ by $y_i$;
then the map $a\otimes \eta_i \mapsto \phi^*a\cup y_i$, $a\in H^\bullet(X/G;\Q)$
extends to an isomorphism of graded $\Q$-vector spaces
$$
H^\bullet (X/G;\Q)\otimes H^\bullet(G;\Q) \iso H^\bullet(X;\Q). 
$$
\end{theorem}

\section{Properties of fundamental classes} \label{fundclasses}

We collect some facts on fundamental classes
that we need later on. We refer to \cite{E} for the cohomology-version and
\cite{F} for the Chow-version. 
\vspace{1ex}
 
\noindent 1. For any connected submanifold $Z$ of pure codimension $c$ in a
complex algebraic manifold $X$, its fundamental class $[Z]\in H^{2c}_Z(X)$ is the
image of $1\in H^0(Z)$ under the Thom-isomorphism $H^\bullet(Z)
\iso~ H^\bullet_Z(X)[2c]$. For $Z$ an irreducible
subvariety, one still has a fundamental class as above, since restriction to
the smooth part of $Z$ induces isomorphisms between the  relevant cohomology
groups with support in $Z$, respectively the smooth part of $Z$. If $Z=\sum_i
n_i Z_i$ is a cycle of codimension $c$ (with $Z_i$ irreducible), with support
$|Z|$, there is a cycle class $[Z]\in H^{2c}_{|Z|}(X)$. More generally still,
one may assume $Z$ to be a complex subscheme of pure codimension $c$ with
irreducible components $Z_i$ of multiplicity $n_i$ in $Z$ and define the
fundamental class to be the fundamental class of the associated cycle $\sum_i
n_iZ_i$. There are natural maps $H^\bullet_{Z_i} \to H^\bullet_{|Z|}$ and if we
identify $[Z_i]$ with their images under these maps we have the  equality
$$
[Z]=\sum_i n_i[Z_i].
$$

\noindent 2. The fundamental classes behave functorially as follows. Let $f:X\to
Y$ be a holomorphic map between complex algebraic manifolds, $Z\subset X$,
$W\subset Y$ subschemes such that $Z$ is contained in the scheme-theoretic
inverse image $f^{-1}W$. Then $f$ induces $H^\bullet_W(Y)\to H^\bullet_Z(X)$
and if moreover $Z=f^{-1}W$ has the same codimension $c$ as $W$, then
$f^*[W]=[Z]$. In particular, if $W$ is irreducible and the cycle
associated to $Z=f^{-1}W$ is $\sum n_iZ_i$, we find
$$
f^*[W]= [f^{-1}W]=\sum n_i [Z_i] \in H^{2c}_{|Z|}(X).
$$

\noindent 3. We can refine the fundamental class of $Z$, a purely
$c$-codimensional subscheme of $X$ to a class in the Chow group $A_{n-c}(X)$,
$n=\dim(X)$. The refinement works as follows. There is a cycle class map
$A_\bullet(X)\to H^{\rm BM}_{2\bullet}(X)$ to Borel-Moore homology. One composes
this map in degree $n-c$ with Poincar\'e-duality for Borel-Moore homology, which
reads
$$
H^{\rm BM}_\bullet(X) \stackrel{\sim}{\rightarrow} H^{2n-\bullet}_Z(X).
$$
In sum, we get a cycle class map
$$
A_{n-c}(X) \to H^{2c}_Z(X)
$$
sending the Chow cycle of $Z$ to $[Z]$. Abusing notation, we denote the Chow
cycle also by $[Z]$. This is especially useful if $Z$ is the scheme of zeros of
a section $s$ of a vector bundle $E$ over $X$. In fact, if $s:E \to X$ is the
zero-section with image, say $\{0\}$, there is a Gysin isomorphism $s^*:
A_\bullet(E) \to A_\bullet(X)[-r]$ with the property
$$
A_n(E) \ni [\{0\}]  \stackrel{\displaystyle s^*}{\longmapsto}  c_r(X)\in
A_{n-r}(X).
$$
See \cite[Example 3.3.2]{F}. This Gysin map is in fact the inverse of the
isomorphism
$$
\pi^* :A_{n-r}(X)\iso A_{n}(E).
$$

\section{The cohomology ring of the general linear group} \label{gln}

We turn to $G=G_n=\mathrm{GL}_n(\C)$, $n\ge 1$. In this case, by \cite{B},
$H^\bullet(G)$ is the exterior algebra with generators $\eta_\ell^{(n)}$in
all odd degrees $2\ell-1$, $\ell=1,\dots,n$. In other words
$r_1=1,r_2=2,\dots,r_n=n$. Since
$G_n
\subset  M_n=\mathrm{Mat}_n(\C)$, a vector space, we have an identification of
mixed Hodge structures
$$H^\bullet(G) \iso H^\bullet_{D_n}(M_n)[1], 
$$
where
$$
D_n = \{A \in M_n \mid \det(A)=0\} = M_n \setminus G_n,
$$
and so $\eta_\ell^{(n)}$ corresponds to some class in $H^{2\ell}_{D_n}(M_n)$.
The goal is to find explicit descriptions of this class as fundamental
class   of the subvariety $D_{n,\ell} \subset D_n$ to be defined below. This 
will turn out to be essential for the next section. We are going to
show this by first defining classes $\eta_\ell^{(n)}$ that clearly have this
property. Then we prove that these classes do generate $H^\bullet(G)$ as an
exterior algebra. 
\vskip 1ex
\noindent We introduce the following notation: 
\begin{itemize}
\item $D_{n,\ell} \subset D_n$: the subvariety consisting of
those matrices for which the first $n+1-\ell$ columns are linearly dependent.
 Note that $D_{n,\ell}$ has codimension $\ell$ in $M_n$. 
\item $\tilde{D}_n = \{(A,p) \in D_n \times \PP^{n-1}(\C) \mid A(p)=0\}$  and
$\pi_n:\tilde{D}_n \to D_n$ is the projection to the first factor. 
\item $Q_n = \{(x,y) \in \C^n\times \C^n \mid x\bullet y = 1\}$.
\item $\alpha_n:M_{n-1}\to M_n$ is the inclusion which maps a matrix $A$ to 
$\left (
\begin{array}{cc}1 & 0 \\ 0 & A \end{array} \right )$. 
\item $h$: the hyperplane class in $H^2(\PP^n(\C))$.  
\end{itemize}
Note that the projection to the second factor turns $\tilde{D}_n$  into a
vector bundle of rank $n^2-n$ over $\PP^{n-1}(\C)$, so $\tilde{D}_n$ is smooth
and $\pi_n$ is a resolution of singularities of $D_n$. 

\begin{lemma}
Let $X$ be a smooth variety, $D\subset X$ a subvariety of  codimension $k$ and
$\pi:\tilde{D} \to D$ a resolution of singularities. Then  there are natural
Gysin maps $\beta_\ell:H^{\ell-2k}(\tilde{D})(-k) \to H^{\ell}_D(X)$ which are
morphisms of mixed Hodge structures. 
\end{lemma}
{\em Proof } Let $n=\dim(X)$.  Because $X$ is smooth, cap product with the 
fundamental class $[X] \in H^{\rm BM}_{2(n-k)}(X)$  in Borel-Moore homology
induces isomorphisms of mixed Hodge structures 
$$
H^{\ell}_D(X) \simeq H^{\rm BM}_{2n-\ell}(D)(-n).
$$
by \cite[Sect.~19.1]{F}. As Borel-Moore homology is covariant for proper morphisms  we have natural
maps
$$ H^{\rm BM}_j(\tilde{D})\to H^{\rm
BM}_j(D).
$$ 
As $\tilde{D}$ is smooth, cup product with the fundamental class $[\tilde{D}]$
induces an isomorphism 
$$
H^{\ell-2k}(\tilde{D}) \to H^{\rm BM}_{2(n-k)-\ell}(\tilde{D})(k-n) .$$ 
The map $\beta_\ell$ is obtained the isomorphism 
$$H^{\ell-2k}(\tilde{D})(-k) \to H^{\rm BM}_{2(n+k)-\ell}(\tilde{D})(-n)$$ foolowed by the direct image map 
$$H^{\rm BM}_{2(n+k)-\ell}(\tilde{D})(-n) \to H^{\rm BM}_{2(n+k)-\ell}(D)(-n)$$
and the inverse of the isomorphism 
$$
H^{\ell-2k}_D(X)(-k) \simeq H^{\rm BM}_{2n+2k-\ell}(D)(-n).
$$
$\square$

\vspace{4mm}
\noindent Let us apply this to the situation of  $\tilde{D}_n \to D_n
\hookrightarrow M_n$.   We obtain maps 
$$
\beta^{(n)}_\ell: H^{2\ell-2}(\PP^{n-1}(\C))(-1) \to
H^{2\ell}_{D_n}(M_n))\simeq H^{2\ell-1}(G_n)
$$
and define for $\ell = 1,\ldots,n$: 
$$
\eta^{(n)}_\ell := \beta^{(n)}_\ell\left(\frac{h^{\ell-1}}{2\pi i}\right)\in
H^{2\ell-1}(G_n).$$ 

We observe that the class in $H_{D_n}^{2\ell}(M_n)$
corresponding to $\eta_\ell^{(n)}$ is indeed the fundamental class of
$D_{n,\ell} \subset D_n$.

\begin{lemma}
The map $\alpha:M_{n-1} \to M_n$ maps  $D_{n-1}$ and $G_{n-1}$ to $D_n$ and
$G_n$  respectively and $\alpha^\ast(\eta^{(n)}_\ell) = \eta^{(n-1)}_\ell$ for
$\ell = 1,\ldots,n-1$ while  $\alpha^\ast(\eta^{(n)}_n)=0$. 
\end{lemma}  {\em Proof } Observe that $\alpha^{-1}(D_{n,\ell}) =
D_{n-1,\ell}$.   One checks that this holds not only set theoretically, but even
as schemes. Then the lemma follows from property 2) from section
\ref{fundclasses}.

Because the classes $\eta^{(n)}_\ell$ are of odd degree, they have square zero
and anti-commute, so we have a homomorphism of graded algebras 
$$
R_n: \Lambda(z_1,\ldots,z_n) \to H^\ast(G_n).
$$  
Here $\Lambda(z_1,\ldots,z_n)$ is the exterior algebra on $n$ generators
$z_1,\ldots,z_n$ with $z_i$ of degree $2i-1$, and
$R_n(z_\ell)=\eta^{(n)}_\ell$. 
\begin{theorem}
The map $R_n$ is an isomorphism. Moreover, the generators $\eta_\ell^{n}\in
H^{2\ell-1}(G_n)$ have pure type $(\ell,\ell)$ and map to the fundamental classes
$D_{n,\ell}$ under the identification $H^{2\ell-1}(G_n)\simeq
H^{2\ell}_{D_n}(M_n)$.
\end{theorem}
{\em Proof } By induction on $n$. For $n=1$ everything is clear. Suppose the map
 $R_{n-1}$ is an isomorphism. We consider the map 
$$
\rho:G_n \to Q_n, \ \rho(g) = (g(e_1), ^tg^{-1}(e_1)).
$$
This is  the orbit map of a transitive action of $G_n$ on $Q_n$ and
$\alpha(G_{n-1})$ is the isotropy subgroup of $(e_1,e_1)\in Q_n$. Therefore,
$\rho$ is also the quotient map for the action of $G_{n-1}$ on $G_n$ by left
translation via $\alpha$. As the classes $\eta^{(n-1)}_\ell$ generate the
cohomology ring of $G_{n-1}$ and are images of classes on $G_n$, the
restriction maps $\alpha^\ast: H^i(G_n) \to H^i(G_{n-1})$ are surjective. Hence
by Theorem \ref{LH} we have an isomorphism 
$$
H^\ast(Q_n) \otimes H^\ast(G_{n-1}) \simeq H^\ast(G_n).
$$
The variety $Q_n$ is homotopy equivalent to a sphere of dimension $2n-1$ (in  fa
ct to its
subvariety consisting of pairs $(x,y)$ with $y=\overline{x}$).   Moreover, a
generator of $H^{2n-1}(Q_n)$ is mapped to a non-zero multiple of $\eta^{(n)}_n$
by the map $\rho^\ast$. This implies the surjectivity and hence bijectivity of
$R_n$. $\square$

\vspace{4mm}\noindent
{\bf Remark } For any Lie group $G$, the map $g \mapsto g^{-1}$ induces 
multiplication by
$-1$ on the Lie algebra, hence on $H^k(G)$ it induces multiplication by
$(-1)^k$. The involution $\sigma: G_n\to G_n$ given by $\sigma(g) = ^tg^{-1}$
has $\sigma^\ast(\eta^{(n)}_n)=(-1)^n\eta^{(n)}_n$. Indeed, if we let
$\sigma:Q_n\to Q_n$ be given by $\sigma(x,y)=(y,x)$ then $\rho$ becomes
equivariant, and it is an easy exercise to see that $\sigma^\ast = (-1)^n$ on
$H^{2n-1}(Q_n)$. We conclude that transposition $\tau$ on $G_n$ induces
$\tau^\ast(\eta^{(n)}_n)=(-1)^{n-1}\eta^{(n)}_n$. As the inclusion $G_{n-1}\to G_n$ commutes
with transposition, we conclude that
$\tau^\ast(\eta^{(n)}_\ell)=(-1)^{\ell-1}\eta^{(n)}_\ell$ for all $\ell
\leq n$.

\section{Moduli of smooth hypersurfaces} \label{moduli}
We let $\Pi_{n,d} = \C[x_0,\ldots,x_n]_d$ denote the vector  space of
homogeneous  polynomials of degree $d$ in $n+1$ variables over $\C$. We let
$$
\Sigma_{n,d} = \{f \in
\Pi_{n,d} \mid f \mbox{ has a critical point outside } 0\}.
$$ 
There exists an irreducible polynomial $\Delta$ in the coefficients of $f\in
\Pi_{n,d}$ such that $f\in\Sigma_{n,d}$ if and only if $\Delta(f)=0$. Moreover,
$\Delta$ is homogeneous of degree $(n+1)(d-1)^n$. 

We let $U_{n,d} = \Pi_{n,d}\setminus \Sigma_{n,d}$. The group
$\mathrm{GL}_{n+1}(\C)$ acts  on $U_{n,d}$. For $d\leq 2$ or $d=3,n=1$ it acts
transitively, but in the remaining cases it acts with finite isotropy groups
and we have a geometric quotient $M_{n,d}$ which is a coarse moduli space for
non-singular projective hypersurfaces of degree $d$ in $\PP^n(\C)$. 
In our situation we fix a particular $f=f_{n,d}\in U_{n,d}$, the 
Fermat  hypersurface:
$$f_{n,d} = x_0^d +\cdots +x_n^d,$$
and the orbit map then extends to a map 
\begin{eqnarray*} 
r_n: M_{n+1} &\to \Pi_{n,d} \cr
A &\mapsto f_{n,d}\circ A.
\end{eqnarray*}
It induces maps for cohomology with supports:
$$
H^{2\ell}_{\Sigma_{n,d}}(\Pi_{n,d})  \stackrel{r_n^\ast}{\longrightarrow}  
H^{2\ell}_{D_{n+1}}(M_{n+1}). 
$$
Define for $\ell=1,\dots,n+1$
$$
\Sigma^{(\ell)}_{n,d} = \{f\in \Pi_{n,d} \mid V(f)^{\mathrm{sing}}\cap 
[e_0,\ldots,e_{n-\ell+1}] \neq \emptyset\}.
$$
Then  $\Sigma^{(\ell)}_{n,d}\subset \Sigma_{n,d}$  has codimension  $\ell$ in
$\Pi_{n,d}$. Below we shall prove:
\begin{lemma} \label{mainlemma} The class $r_n^\ast([\Sigma^{(\ell)}_{n,d}])$ is
a  non-zero  multiple of $[D_{n+1,\ell}]$.
\end{lemma}
Recall from the previous section that $[D_{n+1,\ell}]$ corresponds to the
generator $\eta^{(n)}_\ell\in H^{2\ell -2}(G)$ and we now apply Theorem
\ref{maintheorem} to deduce:

\begin{theorem} \label{leraymoduli}
Let $d \geq 3$. Then the Leray spectral sequence in  rational cohomology for
the quotient  map $U_{n,d} \to M_{n,d}$ degenerates at $E_2$.
\end{theorem}
 
Let us proceed to give a {\bf proof of Lemma} \ref{mainlemma}.
We want to do this by induction on $n$, so we fix an embedding
 $\iota:\Pi_{n-1,d} \hookrightarrow \Pi_{n,d}$ by posing 
$$
\iota(h) = x_0^d+h(x_1,\ldots,x_n).
$$
Note that $\iota(f_{n-1,d})=f_{n,d}$ and that $\iota(\Pi_{n-1,d})\cap 
\Sigma_{n,d}=\iota(\Sigma_{n-1,d})$. The intersection  multiplicity however is equal
to $d-1$. Indeed, the multiplicity of a stratum of the discriminant
corresponding to hypersurfaces with isolated singularities is equal to the sum
of their Milnor numbers, and adding the term $x_0^d$ multiplies the Milnor
numbers by $d-1$. We obtain a commutative diagram 
$$
\begin{array}{ccc}
M_n & \stackrel{r_{n-1}}{\longrightarrow} & \Pi_{n-1,d} \\
\downarrow \alpha & & \downarrow \iota\\
M_{n+1} & \stackrel{r_n}{\longrightarrow} & \Pi_{n,d}. 
\end{array}
$$
We have a corresponding diagram in 
cohomology with supports 
$$
\begin{array}{ccc}
H^{2\ell}_{\Sigma_{n,d}}(\Pi_{n,d}) & \stackrel{r_n^\ast}{\longrightarrow} & 
H^{2\ell}_{D_{n+1}}(M_{n+1}) \\
\downarrow \iota^* & & \downarrow \alpha^* \\
H^{2\ell}_{\Sigma_{n-1,d}}(\Pi_{n-1,d}) & \stackrel{r_n^\ast}{\longrightarrow}
& H^{2\ell}_{D_n}(M_n)
\end{array}
$$
Observe that $\iota^\ast([\Sigma^{(\ell)}_{n,d}]) =
\nu^{(\ell)}_{n,d}[\Sigma^{(\ell)}_{n-1,d}]$ where $\nu^{(\ell)}_{n,d}$ is the 
intersection multiplicity of $\Sigma^{(\ell)}_{n,d}$ with $\iota(\Pi_{n-1,d})$ in
$\Pi_{n,d}$. In particular, $\nu^{(\ell)}_{n,d}$ is a positive integer. 

We can now prove the lemma  by induction on $n$ using the above
diagram, provided we check the case $\ell=n+1$ for each $n$. 

The variety $S=\Sigma^{(n+1)}_{n,d}$ is the  linear space of all polynomials
singular at  $e_0$. Its pre-image under $r_n$  has two irreducible components:
one consists of the  matrices whose first column is zero, i.e. with
$A(e_0)=0$;  this component is exactly $T_1= D_{n+1,n+1}$. The other component,
$T_2$,  which has the same dimension, consists generically of matrices $A$
mapping $e_0$  to some point $p$ with $f_{n,d}(p)=0$ and such that the image of $A$
is contained  in the tangent space to the hypersurface $V(f_{n,d})\subset
\C^{n+1}$ at this point. The   component $T_2$ has  multiplicity one, whereas
$T_1$  has multiplicity  $d(d-1)^n$. We have the commutative diagram 
$$
\begin{array}{ccc}
H^{2n+2}_S(\Pi_{n,d}) & \rightarrow & H^{2n+2}_{\Sigma_{n,d}}(\Pi_{n,d}) \\
\downarrow r_n^\ast & & \downarrow \\
H^{2n+2}_{T_1\cup T_2}(M_{n+1}) & \rightarrow & H^{2n+2}_{D_{n+1}}(M_{n+1})
\end{array}
$$
and therefore 
$$
r_n^\ast([S])=d(d-1)^n[T_1]+[T_2]  \leqno{(1)}  
$$
by Property 1) in Sect.~
\ref{fundclasses}.

\noindent {\bf Claim}: {\it We have
$$[T_2]= (-1)^n(1-(1-d)^n)[T_1] \mbox{ in
}H^{2n+2}_{D_{n+1}}(M_{n+1}).
$$  
}
Combining the Claim with (1) we find:
$$
r_n^\ast[S] = d(d-1)^n[T_1]+[T_2] = ((d-1)^{n+1}+(-1)^n)[T_1] \neq 0,
$$
which proves the Lemma.

It remains to prove the Claim. Let $T_2^\prime$ denote the image of $T_2$ under
the  transposition map $\tau$. Then 
$$
[T_2]= (-1)^n[T_2^\prime] \leqno{(2)}
$$ in
$H^{2n+2}_{D_{n+1}}(M_{n+1})$ by the Remark at the end of Sect.~\ref{gln}. Let
$\tilde{T}_1=T_1\times\{e_0\}\subset \tilde{D}_{n+1}$. 

Write $X=V(f_{n,d}) \subset \PP^n$ and let $\gamma:X \to \PP^n$ be the Gauss  map,
which associates to a point $p \in X$ the coordinates of its tangent
hyperplane, i.e.  $\gamma(p) =
\nabla f_{n,d}(p)$. 

The space $\tilde{D}_{n+1}$ is the total space of a vector bundle $E$ over
$\PP^n$ of  rank $r=n(n+1)$. Let 
$$
\tilde{T} := \{(A,p)\in M_{n+1}\times X \mid  (df_0)_p \circ ^tA=0 \mbox{ and }^tA(e_0)=p\}.
$$
Then $\tilde{T}$ is the total space of a vector bundle $F$ over $X$ of rank
$r-n+1$  which is a subbundle of $\gamma^\ast(E)$, because $(A,p)\in\tilde{T}$
implies that
$A(\gamma(p))=0$. The projection of $\tilde{T}$ in $M_{n+1}$ is precisely
$T_2^\prime$.

We will carry out our calculations in Chow groups instead of cohomology
groups, using property 3) in Sect.~\ref{fundclasses}.  Consider  the
diagram 
$$
\begin{array}{ccccc}
F & \hookrightarrow & \gamma^\ast(E) & \stackrel{\tilde{\gamma}}{\rightarrow} & 
E \\
\downarrow & & \downarrow \pi' & & \downarrow \pi \\
X & = & X & \stackrel{\gamma}{\rightarrow} & \PP^n 
\end{array}
$$
We let $s$ be the $0$-section of $E$, and $s'$ that of $\gamma^*E$ and recall
from Sect.~\ref{fundclasses} 3) that these induce Gysin maps in Chow groups.

The strategy is to compare the classes $\tilde{T_1}$ and $\tilde{T}$ by pushing
them to $\PP^n$. We get two $0$-cycles on $\PP_n$ whose degrees we compare.
Clearly $\deg s^\ast[\tilde{T}_1]=1$ and so it suffices to calculate the degree 
of
$$
s^\ast\tilde{\gamma}_\ast([F]) \in
A_0(\PP^n).
$$
By \cite[Proposition 1.7]{F} we find that 
$$
 \tilde{\gamma}_\ast\pi^{\prime\ast}\alpha=\pi^\ast\gamma_\ast\alpha \in
A_{i+r}(E)
$$
for any $\alpha \in A_i(X)$. Applying  this to $\alpha=s^{\prime\ast}[F]$   we
find 
$$
\tilde{\gamma}_\ast[F] = \pi^\ast\gamma_\ast s^{\prime\ast}[F].
$$ 
Next, applying $s^\ast$ to both sides and using that the Gysin map $s^*$ is in
fact the inverse of the isomorphism induced by the bundle projection $\pi:E\to
\PP^n$, and similarly for $s'$, we get 
$$
s^\ast\tilde{\gamma}_\ast[F]=\gamma_\ast s^{\prime\ast}[F].
$$ 
We next compute $s^{\prime\ast}[F] \in A_0(X)$. By \cite[Example 3.3.2]{F}  
applied to the vector bundle $\gamma^\ast(E)/F$ we get 
$$ 
s^{\prime\ast}[F]  = c_{n-1}(\gamma^\ast(E)/F).
$$
On $\PP^n$ we have the exact sequence 
$$
0 \to E \to \calO^{(n+1)^2} \to \calO(1)^{n+1} \to 0
$$
showing that $c(E) = (1+h)^{-n-1}$ where $h=c_1(\calO(1))$.  As  $\gamma^\ast
\calO(1) = \calO_X(d-1)$ we get 
$$
c(\gamma^\ast E) = (1+(d-1)h_X)^{-n-1} 
$$
where $h_X = c_1(\calO_X(1))$. For the bundle $F$ we have the exact sequences 
$$
0 \to F \to \calO_X^{(n+1)^2} \to Q_X \oplus \calO_X(d-1)^n \to 0$$
$$
0 \to \calO_X(-1) \to \calO_X^{n+1} \to Q_X \to 0
$$
so $Q_X$ is the restriction of the universal quotient bundle to $X$. Hence we find 
$$c(F) = (1+(d-1)h_X)^{-n}c(Q_X)^{-1} = (1+(d-1)h_X)^{-n}(1-h_X)^{-1}
$$
so 
$$c(\gamma^\ast E/F)= (1+(d-1)h_X)^{-1}(1-h_X)^{-1}.$$
We find 
$$c_{n-1}(\gamma^\ast E/F)= \left(\frac{1-(1-d)^n}{d}\right)h_X^{n-1}$$ 
which has degree equal to $1-(1-d)^n$. Combining this with (2), the Claim then
follows.
\hfill$\square$

\end{document}